\documentclass[11pt]{article}
\usepackage{amssymb}
\usepackage{amsmath}
\usepackage{amsthm}
\usepackage{latexsym}
\usepackage{amsfonts}
\usepackage{graphicx}
\usepackage{graphics}
\usepackage[T1]{pbsi}

\newtheorem{thm}{Theorem}[section]   
\newtheorem{lem}{Lemma}
\newtheorem{cor}{Corollary}

\theoremstyle{definition}
\newtheorem*{Proof}{Proof}

\newcommand{\dis}{\displaystyle}
\newcommand{\fa}{\forall}

\newcommand{\el}{\ell}

\newcommand{\ra}{\;\rightarrow\;}

\newcommand{\bi}{\beta}
\newcommand{\ga}{\gamma }

\newcommand{\De} {{\varDelta}}

\newcommand{\vthi}{\vartheta }

\newcommand{\la}{\lambda }
\newcommand{\mi}{\mu }

\newcommand{\ti}{\tau }

\newcommand{\oo}{\omega}

\newcommand{\R}{\mathbb{R}}
\newcommand{\Z}{\mathbb{Z}}

\newcommand{\cf}{{\cal{F}}}

\newcommand{\ct}{{\cal{T}}}

\newcommand{\cm}{{\cal{M}}}

\newcommand{\ld}{\ldots}

\newcommand{\sm}{\smallsetminus}

 \newcommand{\loc}{\mbox{\footnotesize loc}}

\begin{document}

\title{\bf THE GEOMETRY OF THE DYADIC MAXIMAL OPERATOR}
\author{Eleftherios N. Nikolidakis}
\footnotetext{\hspace{-0.5cm}} \footnotetext{\hspace{-0.5cm}E-mail address:
lefteris@math.uoc.gr}
\date{}
\maketitle
\noindent
{\bf Abstract:} We prove a sharp integral inequality which connects the dyadic maximal operator with the Hardy operator. We also give some applications of this inequality.
%
%
\section{Introduction}
\noindent

The dyadic maximal operator on $\R^n$ is defined by
\begin{eqnarray}
\cm_d\phi(x)=\sup\bigg\{\frac{1}{|Q|}\int_Q|\phi(u)|du:x\in Q,Q\subseteq\R^n \;\text{is a dyadic cube}\bigg\}  \label{eq1.1}
\end{eqnarray}
for every $\phi\in L^1_{\loc}(\R^n)$ where the dyadic cubes are those formed by the grids $2^{-N}\Z^n$, for $N=0,1,2,\ld\,.$

As it is well known it satisfies the following weak type (1,1) inequality:
\begin{eqnarray}
|\{x\in\R^n:\cm_d\phi(u)>\la\}|\le\frac{1}{\la}\int_{\{\cm_d\phi>\la\}}
|\phi(u)|du, \label{eq1.2}
\end{eqnarray}
for every $\phi\in L^1(\R^n)$ and every $\la>0$.

(\ref{eq1.2}) easily implies the following $L^p$ inequality
\begin{eqnarray}
\|\cm_d\phi\|_p\le\frac{p}{p-1}\|\phi\|_p.  \label{eq1.3}
\end{eqnarray}
It is easy to see that the weak type inequality (\ref{eq1.2}) is best possible, while (\ref{eq1.3}) is also sharp. (See \cite{1}, \cite{2} for general martingales and \cite{19} for dyadic ones).

An approach for studying the dyadic maximal operator is by making certain refinements of the above inequalities. Concerning (\ref{eq1.2}), some of them have been done in \cite{6}, \cite{10}, \cite{11}, \cite{12}, while for (\ref{eq1.3}) the Bellman function of this operator has been explicetely computed in \cite{3}. It is defined by the following way: For every $f,F,L$ such that $0<f^p\le F$, $L\ge f$ the Bellman function of three variables associated to the dyadic maximal operator is defined by:
\begin{align}
S_p(f,F,L)=\sup\bigg\{&\frac{1}{|Q|}\int_Q(\cm_d\phi)^p:\frac{1}{|Q|}
\int_Q\phi(u)du=f,\nonumber\\
&\frac{1}{|Q|}\int_Q\phi(u)^pdu=F,\;\sup_{R:Q\subseteq R}
\frac{1}{|R|}\int_R\phi(u)du=L\bigg\},  \label{eq1.4}
\end{align}
where $Q$ is a fixed dyadic cube, $R$ runs over all dyadic cubes containing $Q$, and $\phi$ is nonnegative in $L^p(Q)$.

Actually the above calculations have been done in a more general setting. More precisely we define for a non-atomic probability measure space $(X,\mi)$ and a tree $\ct$ the dyadic maximal operator associated to $\ct$ by the following way:
\begin{eqnarray}
\cm_\ct\phi(x)=\sup\bigg\{\frac{1}{\mi(I)}\int_I|\phi|d\mi:x\in I\in\ct\bigg\},  \label{eq1.5}
\end{eqnarray}
for every $\phi\in L^1(X,\mi)$.

In fact, the inequalities (\ref{eq1.2}) and (\ref{eq1.3}) remain true and sharp even in this setting.

Then the respective main Bellman function of two variables is defined by the following way:
\begin{eqnarray}
B_p(f,F)=\sup\bigg\{\int_X(\cm_{\ct}\phi)^pd\mi:\phi\ge0,\;\int_X\phi d\mi=f,\; \int_X\phi^pd\mi=F\bigg\},  \label{eq1.6}
\end{eqnarray}
for $0<f^p\le F$.

It is proved in \cite{3} that (\ref{eq1.6}) equals
\[
B_p(f,F)=F\oo_p(f^p/F)^p, \ \ \text{where} \ \ \oo_p:[0,1]\ra\bigg[1,\frac{p}{p-1}\bigg]
\]
denote the inverse function $H^{-1}_p$ of $H_p$, which is defined by $H_p(z)=-(p-1)z^p+pz^{p-1}$, for $z\in\big[1,\frac{p}{p-1}\big]$. As an immediate result we have that $B_p(f,F)$ is independent of the tree $\ct$ and the measure space $(X,\mi)$.

Actually using this we can compute the following Bellman function of three variables defined by:
\begin{align}B_p(f,F,k)=\sup\bigg\{&\int_K(\cm_{\ct}\phi)^pd\mi:\phi\ge0,\;
\int_X\phi d\mi=f,\;\int_X\phi^pd\mi=F,\nonumber \\
&K\;\text{measurable subset of $X$ with}\;\mi(K)=k\bigg\},  \label{eq1.7}
\end{align}
for $0<f^p\le F$ and $k\in(0,1]$.
 
Using(\ref{eq1.6}) one can also find the exact value of (\ref{eq1.4}).

There are several problems in Harmonic Analysis where Bellman functions arise. Such problems (including the dyadic Carleson imbedding theorem and weighted inequalities) are described in \cite{9} (see also \cite{7}, \cite{8}) and also connections to Stochastic Optimal Control are provided, from which it follows that the corresponding Bellman functions satisfy certain nonlinear second-order PDEs.

The exact evaluation of a Bellman function is a difficult task which is connected with the deeper structure of the corresponding Harmonic Analysis problem.

Until now several Bellman functions have been computed (see \cite{1}, \cite{2}, \cite{3}, \cite{7}, \cite{15}, \cite{16}, \cite{17}, \cite{18}).

Recently L. Slavin,A. Stokolos and V.Vasyunin (\cite{14}) in some cases linked the Bellman function computation to solving certain PDEs of the Monge-Amp\`{e}re type, and in this way they obtained an alternative proof of the results in \cite{3} for the Bellman functions related to the dyadic maximal operator. Also in \cite{18} using the Monge-Amp\`{e}re equation approach a more general Bellman function that the one related to the dyadic Carleson Imbedding Theorem has been precisely evaluated.

Also the Bellman functions of the dyadic maximal operator in relation with Kolmogorov's inequality have been evaluated in \cite{5}.

In \cite{4} now more general Bellman functions have been computed such as:
\begin{align}
T_{p,G}(f,F,k)=\sup\bigg\{&\int_KG(\cm_\ct\phi)d\mi:\phi\ge0,\;\int_X\phi d\mi=f,\;\int_X\phi^pd\mi=F,\nonumber \\
&K\;\text{measurable subset of $X$ with}\;\mi(K)=k\bigg\} \label{eq1.8}
\end{align}
where $G$ is a suitable non-negative increasing convex function on $[0,+\infty)$ . For example one can use  $G(x)=x^q$,with  $1<q<p$.

The approach for evaluating (\ref{eq1.8}) is by proving a symmetrization principle, namely that for suitable $G$ as above the following holds
\begin{align}
T_{p,G}(f,F,k)=\sup\bigg\{&\int^k_0G\bigg(\frac{1}{u}\int^u_0r(t)dt\bigg)du:\;r\ge0,\;
r\;\text{non increasing}\nonumber\\
&\text{on $[0,1]$ and}\; \int^1_0r(u)du=f,\;\int^1_0r^p(u)du=F\bigg\}  \label{eq1.9}
\end{align}
Equation (\ref{eq1.9}) is of much importance and is the tool for finding the exact value of $T_{p,G}(f,F,k)$ as is done in \cite{4}.

In this paper we prove a sharp integral inequality which connects the dyadic operator with the Hardy operator in an immediate way.

In fact we consider non-increasing integrable functions $g,h:(0,1]\ra\R^+$, and a nondecreasing function $G:[0,+\infty)\ra[0,+\infty)$. We prove the following
\begin{thm}\label{thm1.1}
\begin{align}
&\sup\bigg\{\int_KG[(\cm_\ct\phi)^\ast]h(t)dt,\;\phi^\ast=g,\;K\;\text{measurable subset of $[0,1]$ with}\; |K|=k\bigg\} \nonumber\\
&=\int^k_0G\bigg(\frac{1}{t}\int^t_0g(u)du\bigg)h(t)dt, \ \ \text{for any} \ \ k\in(0,1].  \label{eq1.10}
\end{align}
\hfill$\square$
\end{thm}

An immediate consequence of the above theorem is the following
\begin{cor}\label{cor1.1}
With the above notation we have that
\[
\sup\bigg\{\int_X(\cm_\ct\phi)^pd\mi:\phi^\ast=g\bigg\}=\int^1_0\bigg(
\frac{1}{t}\int^t_0g(u)du\bigg)^pdt.
\]
for any $p>0$.
\end{cor}

It is obvious that the above theorem implies the symmetrization principle mentioned above.

We believe that Theorem \ref{thm1.1} has many and important applications in the theory of the dyadic maximal operator. We describe some of them as follows:

First of all it is interesting to see what happens if in (\ref{eq1.8}) we set $G(x)=x^q$ and replace the $L^p$-norm of $\phi$ by its $L^{p,\infty}$-quasi norm $\|\cdot\|_{p,\infty}$ defined by
\begin{eqnarray}
\|\phi\|_{p,\infty}=\sup\{\mi(\{\phi\ge\la\})^{1/p}\cdot\la:\la>0\}.  \label{eq1.11}
\end{eqnarray}
More precisely using Theorem \ref{thm1.1} we can evaluate the following
\begin{align}
\De(f,F,k)=\sup\bigg\{&\int_K(\cm_\ct\phi)^qd\mi:\phi\ge0,\;\int_X\phi d\mi=f,\;\|\phi\|_{p,\infty}=F,\nonumber \\
&K\;\text{measurable subset of $X$ with}\; \mi(K)=k\bigg\},  \label{eq1.12}
\end{align}
for every $0<f\le\frac{p}{p-1}F$, $k\in[0,1]$ and $1<q<p$.

Secondly it is known by \cite{10} that the following inequality
\begin{eqnarray}
\|\cm_\ct\phi\|_{p,\infty}\le\frac{p}{p-1}\|\phi\|_{p,\infty},  \label{eq1.13}
\end{eqnarray}
has been proved to be best possible and independent of the $L^1$ and $L^q$-norm of $\phi$, for any fixed $q$ such that $1<q<p$. In \cite{20} it is introduced a norm $L^{p,\infty}$ equivalent to $\|\cdot\|_{p,\infty}$. This is given by
\begin{align}
\||\phi\||_{p,\infty}=\sup\bigg\{\mi(E)^{-1+\frac{1}{p}}\int_E|\phi|d\mi:\;&E
\;\text{measurable subset of}\nonumber\\
 &X\;\text{with}\;\mi(E)>0\bigg\}  \label{eq1.14}
\end{align}
and it is easily proved that the following holds:
\begin{eqnarray}
\|\phi\|_{p,\infty}\le\||\phi\||_{p,\infty}\le\frac{p}{p-1}\|\phi\|_{p,\infty}.
\label{eq1.15}
\end{eqnarray}
As a second application we prove that the following inequality:
\begin{eqnarray}
\||\cm_\ct\phi\||_{p,\infty}\le\bigg(\frac{p}{p-1}\bigg)^2\|\phi\|_{p,\infty},
\label{eq1.16}
\end{eqnarray}
is best possible and independent of the $L^1$-norm of $\phi$. At last we prove that the inequality  $\|\cm_\ct\phi\|_{L^{p,q}}\le\frac{p}{p-1}\|\phi\|_{L^{p,q}}$ is best possible
for $q<p$ where $\|\cdot\|_{L^{p,q}}$ stands for the Lorentz quasi norm on $L^{p,q}$ given by
\begin{eqnarray}
\|\phi\|_{L^{p,q}}\equiv\|\phi\|_{p,q}=\bigg(\int^1_0[\phi^\ast(t)t^{1/p}]^q
\frac{dt}{t}\bigg)^{1/q}.  \label{eq1.17}
\end{eqnarray}
\section{Preliminaries} 
\noindent

Let $(X,\mi)$ be a non-atomic probability measure space.
\begin{def}\label{def2.1}
A set $\ct$ of measurable subsets of $X$ will be called a tree if the following conditions are satisfied:
\begin{enumerate}
\item[1.] $X\in\ct$ and for every $I\in\ct$ we have that $\mi(I)>0$.
    \item[2.] For every $I\in\ct$ there corresponds a finite or countable subset $C(I)\subseteq\ct$ containing at least two elements such that
\begin{itemize}
\item[(a)] the elements of $C(I)$ are disjoint subsets of $I$
\item[(b)] $I=\cup C(I)$.
\end{itemize}
\item[3.] $\ct=\dis\bigcup_{m\ge0}\ct_{(m)}$ where $\ct_{(0)}=\{X\}$ and $\ct_{(m+1)}=\dis\bigcup_{I\in\ct_{(m)}}C(I)$.
\item[4.] We have that
\[
\lim_{m\ra\infty}\sup_{I\in\ct_{(m)}}\mi(I)=0.
\]
\hfill$\square$

\end{enumerate}
\end{def}

Examples of trees are given in \cite{3}.

The most known is the one given by the family of all dyadic subcubes of $[0,1]^m$.

The following has been proved in \cite{3}.
\begin{lem}\label{lem2.1}
For every $I\in\ct$ and every a such that $0<a<1$ there exists a subfamily $\cf(I)\subseteq\ct$ consisting of disjoint subsets of $I$ such that
\[
\mi\bigg(\bigcup_{J\in\cf(I)}J\bigg)=\sum_{J\in\cf(I)}\mi(J)=(1-a)\mi(I).
\]
\hfill$\square$
\end{lem}
We will need also the following fact
\begin{lem}\label{lem2.2}
Let $\phi:(X,\mi)\ra\R^+$ and $(A_j)_j$ a measurable partition of $X$ such that $\mi(A_j)>0$ $\fa\, j$. Then if $\int_X\phi d\mi=f$ there exists a rearrangement of $\phi$, say $h$ $(h^\ast=\phi^\ast)$ such that $\frac{1}{\mi(A_j)}\int_{A_j}hd\mi=f$, for every $j$.
\end{lem}
\begin{Proof}
We set $\phi^\ast=g:[0,1]\ra\R^+$.

We find first a measurable set $B_1\subseteq[0,1]$ such that
\setcounter{equation}{0}
\begin{eqnarray}
|B_1|=\mi(A_1) \ \ \text{and} \ \ \frac{1}{|B_1|}\int_{B_1}g(u)du=f. \label{eq2.1}
\end{eqnarray}
Obviously
\begin{eqnarray}
\frac{1}{\mi(A_1)}\int^{\mi(A_1)}_0g(u)du\ge f\ge\frac{1}{\mi(A_1)}
\int^1_{1-\mi(A_1)}g(u)du.  \label{eq2.2}
\end{eqnarray}
As a result there exists $r$ such that $0<r$, $r+\mi(A_1)<1$ and $\frac{1}{\mi(A_1)}\int^{r+\mi(A_1)}_rg(u)du=f$.Then we just need to set  $B_1=[r,r+\mi(A_1)]$.

Then (\ref{eq2.1}) is obviously satisfied.

We define now $h_1:A_1\ra\R^+$ such that $(h_1)^\ast=(g/B_1)^\ast$ which is a function defined on $(0,\mi(A_1))$. Then it is obvious that $\frac{1}{\mi(A_1)}\int_{A_1}h_1=f$. We then continue in the same way for the space $X\sm A_1$ and inductively complete the proof of Lem-\linebreak ma \ref{lem2.2}.  \hfill$\square$
\end{Proof}

Now given a tree $\ct$ on $(X,\mi)$ we define the associated dyadic maximal operator as follows
\[
\cm_\ct\phi(x)=\sup\bigg\{\frac{1}{\mi(T)}\int_I|\phi|d\mi:x\in I\in\ct\bigg\}.
\]
\section{Main Theorem} 
\noindent

Suppose we are given $g,h:(0,1]\ra\R^+$ non increasing integrable functions. Let also $G:[0,+\infty)\ra[0,+\infty)$ be a non decreasing function. We state the following
\setcounter{lem}{0}
\begin{lem}\label{lem3.1}
Let $k\in(0,1]$ and $K$ measurable subset of $(0,1]$ such that $|K|=k$.Then under the above notation the following holds
\[
\int_KG[(\cm_\ct\phi)^\ast]h(t)dt\le\int^k_0G\bigg(\frac{1}{t}\int^t_0g(u)du\bigg)h(t)dt
\]
for every $\phi\in L^1(X,\mi)$ such that $\phi^\ast=g$.
\end{lem}
\begin{Proof}
Let $v$ be the Borel measure on $(0,1]$ defined by $v(A)=\int_Ah(t)dt$, for every $A$ Borel $\subseteq(0,1]$, and set $I=\int_{K}G[(\cm_\ct\phi)^\ast]dv(t)$. Then
\[
I=\int^{+\infty}_{\la=0}v(\{t\in K:(\cm_\ct\phi)^\ast(t)\ge\la\})dG(\la).
\]
Let also $f=\int_X\phi d\mi$. For $0<\la\le f$ we obviously have 
\[
v(\{t\in K:(\cm_\ct\phi)^\ast(t)\ge\la\})=v(K), \ \ \text{since} \ \ (\cm_\ct\phi)^\ast(t)\ge f, \ \ \fa\;t\in[0,1].
\]
Then $I=II+III$, where $II=v(K)[G(f)-G(0)]$ and
\[
III=\int^{+\infty}_{\la=f}v(\{t\in K:(\cm_\ct\phi)^\ast(t)\ge\la\})dG(\la).
\]
Obviously $II\le[G(f)-G(0)]\int^{k}_{0}h(u)du$

Additionally $v(\{t\in K:(\cm_\ct\phi)^\ast(t)\ge\la\})\le v(\{t\in (0,k]:(\cm_\ct\phi)^\ast(t)\ge\la\})$ since $h,(\cm_\ct\phi)^\ast$ are 
nonincreasing and $|K|=k$.

As a consequence $III\le \int^{+\infty}_{\la=f}v(\{t\in (0,k]:(\cm_\ct\phi)^\ast(t)\ge\la\})dG(\la)$
 
Fix now $\la>f$ and let $E_\la=\{\cm_\ct\phi\ge\la\}$. Then there exists a pairwise disjoint family of elements of $\ct$, $(I_j)_j$, such that
\setcounter{equation}{0}
\begin{eqnarray}
\frac{1}{\mi(I_j)}\int_{I_j}\phi d\mi\ge\la, \ \ \text{and} \ \ E_\la=\cup I_j.  \label{eq3.1}
\end{eqnarray}
In fact we just need to consider the family $(I_j)_j$ of elements of $\ct$ maximal under the integral condition (\ref{eq3.1}). From (\ref{eq3.1}) we have that $\int_{I_j}\phi d\mi\ge\la\mi(I_j)$, for every $j$. Since $(I_j)_j$ is pairwise disjoint we have that
\begin{eqnarray}
\int_{E_\la}\phi d\mi\ge\la\mi(E_\la) \ \ \text{so} \ \  \frac{1}{\mi(E_\la)}\int_{E_\la}\phi d\mi\ge\la.  \label{eq3.2}
\end{eqnarray}
Certainly $\int^{\mi(E_\la)}_0\phi^\ast(u)du\ge\int_{E_\la}\phi d\mi$, so (\ref{eq3.2}) gives
\begin{eqnarray}
\frac{1}{\mi(E_\la)}\int^{\mi(E_\la)}_0\phi^\ast(u)du\ge\la.  \label{eq3.3}
\end{eqnarray}
Let now $a(\la)$ be the unique real number on $[0,1]$ such that $\frac{1}{a(\la)}\int^{a(\la)}_0\phi^\ast(u)du=\la$. It's existence is guaranteed by the fact that $\la>f=\int^1_0\phi^\ast(u)du$ (In fact we can suppose without loss of generality that $g(0+)=+\infty$, otherwise we work on $\la\in(f,\|g\|_\infty]$. Notice that if $\|g\|_\infty=A$ and $\phi^\ast=g$, then $\cm_\ct\phi\le A$ a.e. on $X$).

Let also $A_{\la}=\{t\in (0,k]:(\cm_\ct\phi)^\ast(t)\ge\la\}$.

Additionally 
$A_{\la}\subset\{t\in (0,1]:(\cm_\ct\phi)^\ast(t)\ge\la\}=:B_{\la}$,so
$|A_{\la}|\le|B_{\la}|=\mi(E_\la)$.

Let also $\bi(\la)$ be the unique $\bi\in(0,1]$ for which the following holds:
$(0,\bi)\subset A_{\la}$ and such that for every $t>\bi$ we have either
$(\cm_\ct\phi)^\ast(t)<\la$ or $t>k$.
So $A_{\la}$ differs from $(0,\bi)$ except possibly from the endpoint $\bi$.
As a consequence $A_{\la}\subset (0,\bi(\la)]$ and $|A_{\la}|=\bi(\la)$.
From (\ref{eq3.3}) now and the definition of $a(\la)$ we have that 
\[
\frac{1}{\mi(E_\la)}\int^{\mi(E_\la)}_0\phi^\ast(u)du\ge\la=\frac{1}{a(\la)}
\int^{a(\la)}_0\phi^\ast(u)du,
\]
Since $\phi^\ast=g$ is nonincreasing we obtain that $\mi(E_\la)\le a(\la)$.
As a result $|A_{\la}|\le a(\la)$.So $\bi(\la)\le a(\la)$ and consequently 
we have that $A_{\la}\subset (0,a(\la)]$.

But of course  $A_{\la}\subset (0,k]$.

Consequently $A_{\la}\subset\{t\in(0,k]:t\in(0,a(\la)] =\{t\in(0,k]:
\frac{1}{t}\int^t_0g(u)du\ge\la\}$ from the definition of $a(\la)$.

Obviously then
\begin{eqnarray}
III\le\int^{+\infty}_{\la=f}v(\{t\in(0,k]:\frac{1}{t}\int^t_0g(u)du\ge\la\})dG(\la).  \label{eq3.4}
\end{eqnarray}
From the above estimates of II and III we obtain
\begin{eqnarray}
I\le\int^{+\infty}_{\la=0}v(\{t\in(0,k]:\frac{1}{t}\int^t_0g(u)du\ge\la\})dG(\la)=
\int^k_0G\bigg(\frac{1}{t}\int^t_0g(u)du\bigg)dv(t). 
\end{eqnarray}
and Lemma \ref{lem3.1} is proved.  \hfill$\square$
\end{Proof}

We now complete the proof of Theorem \ref{thm1.1}
\begin{thm}\label{thm3.1}
For any $k\in(0,1]$
\begin{align}
&\sup\bigg\{\int_KG[(\cm_\ct\phi)^\ast]h(t)dt,\;\phi^\ast=g,\;K\;\text{measurable subset of $[0,1]$ with}\; |K|=k\bigg\} \nonumber\\
&=\int^k_0G\bigg(\frac{1}{t}\int^t_0g(u)du\bigg)h(t)dt. 
\end{align}
\end{thm}
\begin{Proof}
Because of Lemma \ref{lem3.1} we need only to construct for every $a\in(0,1)$ a $\mi$-measurable function $\phi_a:X\ra\R^+$ such that $\phi^\ast_a=g$ and
\[
\underset{a\ra0^+}{\lim\sup}\int^k_0G[(\cm_\ct\phi_a)^\ast]dv\ge
\int^k_0G\bigg(\frac{1}{t}\int^t_0g(u)du\bigg)dv(t).
\]
We proceed to this as follows:

Let $a\in(0,1)$. Using Lemma \ref{lem2.1} we choose for every $I\in\ct$ a family $\cf(I)\subseteq\ct$ of disjoint subsets of $I$ such that
\begin{eqnarray}
\sum_{J\in\cf(I)}\mi(I)=(1-a)\mi(I).  \label{eq3.5}
\end{eqnarray}
We define $S=S_a$ to be the smallest subset of $\ct$ such that $X\in S$ and for every $I\in S$, $\cf(I)\subseteq S$. We write for $I\in S$, $A_I=I\sm\dis\bigcup_{J\in\cf(I)}J$. Then if $a_I=\mi(A_I)$ we have because of (\ref{eq3.5}) that $a_I=a\mi(I)$. It is also clear that
\[
S=\bigcup_{m\ge0}S_{(m)}, \ \ \text{where} \ \ S_{(0)}=\{X\}, \ \ S_{(m+1)}=\bigcup_{I\in S_{(m)}}\cf(I).
\]
We define also for $I\in S$, rank$(I)=r(I)$ to be the unique integer $m$ such that $I\in S_{(m)}$.

Additionally we define for every $I\in S$ with $r(I)=m$
\begin{eqnarray}
\ga(I)=\ga_m=\frac{1}{a(1-a)^m}\int^{(1-a)^m}_{(1-a)^{m+1}}g(u)du.  \label{eq3.6}
\end{eqnarray}
We also set for $I\in S$
\[
b_m(I)=\sum_{S\ni J\subseteq I\atop r(J)=r(I)+m}\mi(J).
\]
We easily then see inductively that
\begin{eqnarray}
b_m(I)=(1-a)^m\mi(I). \label{eq3.7}
\end{eqnarray}
It is also clear that for every $I\in S$
\begin{eqnarray}
I=\bigcup_{S\ni J\subseteq I} A_J.  \label{eq3.8}
\end{eqnarray}
At last we define for every $m$ the measurable subset of $X$, 
$S_m:=\bigcup_{I\in S_{(m)}}I$.
Now for every $m\ge0$, we choose $\ti^{(m)}_a:S_{m}\sm S_{m+1}\ra\R^+$ such that
\begin{eqnarray}
[\ti^{(m)}_a]^\ast=\big(g/[(1-a)^{m+1},(1-a)^m)\big)^\ast,  \label{eq3.9}
\end{eqnarray}
This is possible since $\mi(S_m\sm S_{m+1})=\mi(S_{m})-\mi(S_{m+1})=b_m(X)-b_{m+1}(X)=(1-a)^m-(1-a)^{m+1}=a(1-a)^m$ and $X$ is non atomic.

We then set $\ti_a:X\ra\R^+$ by $\ti_a(x)=\ti^{(m)}_a(x)$, for $x\in S_{m}\sm S_{m+1}$, so because of (\ref{eq3.9}), $[\ti^{(m)}_a]^\ast=g$.

It is obvious now that $S_{m}\sm S_{m+1}=\dis\bigcup_{I\in S_{(m)}} A_I$ and that
\begin{align}
&\int_{S_{m}\sm S_{m+1}}\ti^{(m)}_ad\mi=\int^{(1-a)^m}_{(1-a)^{m+1}}g(u)du \nonumber \\
&\Rightarrow
\frac{1}{\mi(S_{m}\sm S_{m+1})}\int_{S_{m}\sm S_{m+1}}
\ti_ad\mi=\ga_m.  \label{eq3.10}
\end{align}
Using now Lemma \ref{lem2.2} we see that there exists a rearrangement of $\ti_a/S_{(m)}\sm S_{(m+1)}=\ti^{(m)}_a$, called $\phi^{(m)}_a$ for which $\frac{1}{a_I}\int_{A_I}\phi^{(m)}_a=\ga_m$, for every $I\in S_m$. Define now $\phi_a:X\ra\R^+$ by $\phi_a(x)=\phi^{(m)}_a(x)$, for $x\in S_{(m)}\sm S_{(m+1)}$. Of course $\phi^\ast_a=g$. Let now $I\in S_{(m)}$. Then
\begin{align}
Av_I(\phi_a)&=\frac{1}{\mi(I)}\int_I\phi_ad\mi=\frac{1}{\mi(I)}\sum_{S\ni J\subseteq I}\int_{A_J}\phi_ad\mi\nonumber \\
&=\frac{1}{\mi(I)}\sum_{\el\ge0}\sum_{S\ni J\subseteq I\atop r(J)=r(I)+\el}\int_{A_J}\phi_ad\mi\nonumber \\
&=\frac{1}{\mi(I)}\sum_{\el\ge0}\sum_{S\ni J\subseteq I\atop r(J)=m+\el}\ga_{m+\el}a_J \nonumber \\
&=\frac{1}{\mi(I)}\sum_{\el\ge0}\sum_{S\ni J\subseteq I\atop r(J)=m+\el}a\mi(J)\frac{1}{a(1-a)^{m+\el}}\int^{(1-a)^{m+\el}}_{(1-a)^{m+\el+1}}g(u)du \nonumber\\
&=\frac{1}{\mi(I)}\sum_{\el\ge0}\frac{1}{(1-a)^{m+\el}}\int^{(1-a)^{m+\el}}_{(1-a)^{m+\el+1}}
g(u)du\sum_{S\ni J\subseteq I\atop r(J)=m+\el}\mi(J) \nonumber \\
&=\frac{1}{\mi(I)}\sum_{\el\ge0}\frac{1}{(1-a)^{m+\el}}\int^{(1-a)^{m+\el}}_{(1-a)^{m+\el+1}}
g(u)du\cdot b_\el(I) \nonumber \\
&\overset{\text{(3.7)}}{=}\frac{1}{(1-a)^m}\sum_{\el\ge0}
\int^{(1-a)^{m+\el}}_{(1-a)^{m+\el+1}}g(u)du=\frac{1}{(1-a)^m}
\int^{(1-a)^m}_0g(u)du.  \label{eq3.11}
\end{align}
Now for $x\in S_m\sm S_{m+1}$, there exists $I\in S_{(m)}$ such that $x\in I$ so
\begin{eqnarray}
\cm_\ct(\phi_a)(x)\ge Av_I(\phi_a)=\frac{1}{(1-a)^m}\int^{(1-a)^m}_0g(u)du=:\vthi_m.  \label{eq3.12}
\end{eqnarray}
Since $\mi(S_m)=(1-a)^m$, for every $m\ge0$ we easily see from the above that we have
\[
(\cm_\ct\phi_a)^\ast(t)\ge\vthi_m, \ \ \text{for every} \ \ t\in[(1-a)^{m+1},(1-a)^m), 
\] 
For any $a\in(0,1)$ we now choose $m=m_{a}$ such that 
$(1-a)^{m+1}\le k<(1-a)^{m}$.

So we have $\underset{a\ra0^+}{\lim}(1-a)^{m_a}=k$

We consider two cases:

(A)
\[
\underset{a\ra0^+}{\lim\sup}\int^k_0G[(\cm_\ct\phi_a)^\ast]dv(t)=+\infty
\]

Then Theorem \ref{thm3.1} is obvious , according to Lemma \ref{lem3.1}.

(B)
\[
\underset{a\ra0^+}{\lim\sup}\int^k_0G[(\cm_\ct\phi_a)^\ast]dv(t)<+\infty
\]

Then
\begin{align}
&\int^{(1-a)^{m_a}}_0G[(\cm_\ct\phi_a)^\ast]dv\ge\sum_{l\ge0}\int^{(1-a)^{m_a+l}}_{(1-a)^{m_a+l+1}}
G(\vthi_m)dv \nonumber\\
&=\sum_{l\ge0}G\bigg(\frac{1}{(1-a)^{m_a+l}}\int^{(1-a)^{m_a+l}}_0g(u)du\bigg)
v([(1-a)^{m_a+l+1},(1-a)^{m_a+l})),  \label{eq3.13}
\end{align}

Since now $\underset{a\ra0^+}{\lim}(1-a)^{m_a}=k$ and the right hand side of 
(\ref{eq3.13}) expresses a Riemman sum for the integral
$\int^{(1-a)^{m_a}}_0G[\frac{1}{t}\int^t_0g(u)du]dv(t)$ ,we conclude
because of the monotonicity of $G,\frac{1}{t}\int^t_0g(u)du$ and $h$ that it 
converges to $\int^{k}_0G(\frac{1}{t}\int^t_0g(u)du)dv(t)$.

So that
\[
\underset{a\ra0^+}{\lim\sup}\int^{(1-a)^{m_a}}_0G[(\cm_\ct\phi_a)^\ast]dv\ge
\int^{k}_0G(\frac{1}{t}\int^t_0g(u)du)dv(t)
\]

Further
\[
\int^{(1-a)^{m_a}}_kG[(\cm_\ct\phi_a)^\ast]dv\le
\big(\int^{(1-a)^{m_a}}_kh(u)du\big)G[(\cm_\ct\phi_a)^\ast(k)]
\]

But
\[
\underset{a>0}\sup G[(\cm_\ct\phi_a)^\ast(k)]<+\infty
\]

otherwise
\[
\underset{a\ra0^+}{\lim\sup} G[(\cm_\ct\phi_a)^\ast(k)]=+\infty
\]

which in turn gives
\[
\underset{a\ra0^+}{\lim\sup}\int^k_0G[(\cm_\ct\phi_a)^\ast(t)]dv(t)=+\infty
\]

that is not our case.

As a result
\[
\underset{a\ra0^+}{\lim}\int^{(1-a)^{m_a}}_kG[(\cm_\ct\phi_a)^\ast(t)]dv(t)=0
\]

Theorem \ref{thm3.1} is now proved.  \hfill$\square$
\end{Proof}

We have now the following
\begin{cor}\label{cor3.1}
For any $p>0$ and $g:(0,1]\ra\R^+$ non increasing we have that
\[
\sup\bigg\{\int_X(\cm_\ct\phi)^pd\mi:\phi^\ast=g\bigg\}=\int^1_0\bigg(\frac{1}{t}
\int^t_0g(u)du\bigg)^pdt.
\]
\hfill$\square$
\end{cor}
\begin{Proof}
Obvious since for any $\phi:(X,\mi)\ra\R^+$
\[
\int_X(\cm_\ct\phi)^pd\mi=\int^1_0[(\cm_\ct\phi)^\ast]^pdt.
\]
\hfill$\square$
\end{Proof}

We give now some applications.
\section{Applications} 
(a) {\bf First application:}

We search for
\setcounter{equation}{0}
\begin{align}
\De(f,F,k)=\sup\bigg\{&\int_K(\cm_\ct\phi)^qd\mi:\phi\ge0,\;\int_X\phi d\mi=f,\;\|\phi\|_{p,\infty}=F,\;K \nonumber\\
&\text{measurable}\;\subseteq X\;\text{with}\;\mi(K)=k\bigg\}  \label{eq4.1}
\end{align}
for $0<f\le\frac{p}{p-1}F$ and $1<q<p$.

We prove
\begin{thm}\label{thm4.1}
For $F=\frac{p-1}{p}$ we have
\begin{eqnarray}
\De(f,F,k)=\left\{\begin{array}{l}
                    \frac{p}{p-q}k^{1-\frac{q}{p}}, \ \ k\le f^{p/p-1} \\ [0.5ex]
                    \frac{q(p-1)}{(p-q)(q-1)}f^{p-q/p-1}-\frac{1}{q-1}
k^{1-q}f^q, \ \ f^{p/p-1}\le k\le1,  \label{eq4.2}
                  \end{array}\right.
\end{eqnarray}
for $0<f\le1$.
\end{thm}
\begin{Proof}

Let $\phi$ be as in (\ref{eq4.1}), and $K$ measurable $\subseteq X$ with $\mi(K)=k$. Using Lemma \ref{lem3.1} we have that
\[
\int_K(\cm_\ct\phi)^qd\mi\le\int^k_0\bigg(\frac{1}{t}\int^t_0\phi^\ast(u)du
\bigg)^qdt.
\]
Since $\|\phi\|_{p,\infty}=\frac{p-1}{p}$ we have that $\phi^\ast(u)\le\frac{p-1}{p}u^{-1/p}$, $u\in(0,1]$. So for every $t$ such that $0<t\le k$
\[
\frac{1}{t}\int^t_0\phi^\ast(u)du\le\frac{1}{t}\int^t_0\frac{p-1}{p}u^{-1/p}=t^{-1/p} \ \ \text{and} \ \ \frac{1}{t}\int^t_0\phi^\ast(u)du\le\frac{f}{t}.
\]
Thus, if we set $A(t)=\frac{1}{t}\int^t_0\phi^\ast(u)du$ we have $A(t)\le\min\big\{\frac{f}{t},t^{-1/p}\big\}$, $\fa\;t\in(0,k]$.

Thus, if $k\le f^{p/p-1}$: $\int^k_0[A(t)]^qdt\le\int^k_0t^{-q/p}dt=\frac{p}{p-q}k^{1-\frac{q}{p}}$ while for $f^{p/p-1}<k\le 1$
\begin{align*}
&\int^k_0[A(t)]^qdt\le\int^{f^{p/p-1}}_0t^{-q/p}dt+\int^k_{f^{p/p-1}}
\frac{f^q}{t^q}dt\\
&=\frac{p}{p-q}f^{p-q/p-1}-\frac{1}{q-1}f^qk^{1-q}+
\frac{1}{q-1}f^{q+\frac{p(1-q)}{p-1}}\\
&=\frac{q(p-1)}{(p-q)(q-1)}f^{p-q/p-1}-\frac{1}{q-1}f^qk^{1-q}.
\end{align*}
So we have proved that $\De\big(f,\frac{p-1}{p},K)\le \ct(f,k)$, where $T(f,k)$ is the right side of (\ref{eq4.2}).

We now prove the reverse inequality.

Obviously, we have that
\begin{eqnarray}
\De\bigg(f,\frac{p-1}{p},k\bigg)\ge\int^k_0\bigg(\frac{1}{t}\int^t_0\psi(u)du\bigg)^qdt, \label{eq4.3}
\end{eqnarray}
where $\psi:(0,1]\ra\R^+$ is defined by $\psi(u)=\left\{\begin{array}{cc}
                                                           \frac{p-1}{p}u^{-1/p}, & 0<u\le f^{p/p-1} \\
                                                           0, & f^{p/p-1}<u\le1
                                                         \end{array}\right.$.
Since $\int^1_0\psi(u)du=f$ and $\|\psi\|^{[0,1]}_{p,\infty}=\frac{p-1}{p}$, (\ref{eq4.3}) is obvious because of Theorem \ref{thm3.2}.

But if $\psi$ is as above we have that
\[
\begin{array}{l}
  \dfrac{1}{t}\int^t_0\psi(u)du=\dfrac{f}{t}, \ \ \text{for} \ \ f^{p/p-1}<t\le1 \ \ \text{and} \\ [2ex]
   \dfrac{1}{t}\int^t_0\psi(u)du=t^{-1/p},  \ \ \text{for} \ \ 0<t\le f^{p/p-1}.
\end{array}
\]
From the above calculations we conclude
\[
\De\bigg(f,\frac{p-1}{p},k\bigg)=T(f,k)
\]
and Theorem \ref{thm4.1} is proved. \hfill$\square$
\end{Proof}
\noindent
(b) {\bf Second application:}

In \cite{10} we have proved that
\begin{eqnarray}
\sup\bigg\{\|\cm_\ct\phi\|_{p,\infty}:\phi\ge0,\;\int_X\phi d\mi=f,\;\|\phi\|_{p,\infty}=F\bigg\}=\frac{p}{p-1}F,  \label{eq4.4}
\end{eqnarray}
for $0<f\le\frac{p}{p-1}F$ , that is the inequality $\|\cm_\ct\phi\|_{p,\infty}\le\frac{p}{p-1}\|\phi\|_{p,\infty}$ is sharp and independent of the integral of $\phi$.

A related problem is to find
\[
E(f,F)=\sup\bigg\{\||\cm_\ct\phi\||_{p,\infty}:\phi\ge0,\;\int_X\phi d\mi=f,\;\|\phi\|_{p,\infty}=F\bigg\}
\]
where is the known integral norm $\||\cdot\||_{p,\infty}$ given by (\ref{eq1.14}). In fact we prove
\begin{thm}\label{thm4.2}
With the above notation we have
\begin{eqnarray}
E(f,F)=\bigg(\frac{p}{p-1}\bigg)^2F.  \label{eq4.5}
\end{eqnarray}
\end{thm}
\begin{Proof}
We prove it for $F=\frac{p-1}{p}$. It is obvious that for every $\phi\in L^{p,\infty}$
\[
\||\cm_\ct\phi\||_{p,\infty}\le\bigg(\frac{p}{p-1}\bigg)^2\|\phi\|_{p,\infty}.
\]
Indeed because of (\ref{eq1.15}) and (\ref{eq4.4})
\begin{eqnarray}
\||\cm_\ct\phi\||_{p,\infty}\le\frac{p}{p-1}\|\cm_\ct\phi\|_{p,\infty}\le
\bigg(\frac{p}{p-1}\bigg)^2\|\phi\|_{p,\infty}, \ \ \text{for every} \ \ \phi\in L^{p,\infty}.  \label{eq4.6}
\end{eqnarray}
We prove now that (\ref{eq4.6}) is best possible and independent of the integral of $\phi$.

Let $0<f\le1$. Choose $k_0$ such that $0<k_0\le f^{p/p-1}$. Set
\[
\psi(u):=\left\{\begin{array}{cc}
                  \frac{p-1}{p}u^{-1/p}, & 0<u\le f^{p/p-1} \\
                  0, & f^{p/p-1}<u\le1.
                \end{array}\right.
\]
Then obviously
\begin{align*}
&E\bigg(f,\frac{p-1}{p}\bigg)\ge\sup\bigg\{k^{-1+\frac{1}{p}}_0\int_E
(\cm_\ct\phi)d\mi:E\;\text{measurable}\;\subseteq X\;\text{with}\;\mi(E)=k_0,\;\phi^\ast=\psi\bigg\} \\
&=k^{-1+\frac{1}{p}}_0\int^{k_0}_0\bigg(\frac{1}{t}\int^t_0\psi(u)du\bigg)dt=
\frac{p}{p-1},
\end{align*}
and Theorem \ref{thm4.2} is proved.  \hfill$\square$
\end{Proof}
\noindent
(c) {\bf Third application:}

We give the last application. We know that the Lorentz space $L^{p,q}(X,\mi)\equiv L^{p,q}$ is defined as
\[
L^{p,q}=\bigg\{\phi:(X,\mi)\ra\R^+ \ \ \text{such that} \ \ \int^1_0[\phi^\ast(t)t^{1/p}]^q\frac{dt}{t}<+\infty\bigg\}
\]
with topology endowed by the quasi-norm $\|\cdot\|_{p,q}$ given by
\[
\|\phi\|_{p,q}=\bigg[\int^1_0[\phi^\ast(t)t^{1/p}]^q\frac{dt}{t}\bigg]^{1/p}.
\]

We prove now the following
\begin{thm}\label{thm4.3}
$\cm_\ct$ maps $L^{p,q}$ to $L^{p,q}$ and $\|\cm_\ct\|_{L^{p,q}\ra L^{p,q}}=\frac{p}{p-1}$, where $q<p$.
\end{thm}
\begin{Proof}
We set $v(A)=\int_Ah(t)dt$, for all Borel subsets $A$ of $[0,1]$, where $h(t)=t^{\frac{q}{p}-1}$. Then
\begin{align}
\|\cm_\ct\phi\|^q_{p,q}&=\int^1_0[\cm_\ct\phi)^\ast t^{1/p}]^q\frac{dt}{t}=\int^1_0[(\cm_\ct\phi)^\ast]^qdv(t) \nonumber\\
&\le\int^1_0\bigg(\frac{1}{t}\int^t_0\phi^\ast(u)du\bigg)^qdv(t).  \label{eq4.7}
\end{align}
We set $A(t)=\frac{1}{t}\int^t_0\phi^\ast(u)du$. Then $A(t)=\int^1_0\phi^\ast(tu)du$. So by the continuous form of Minkowski inequality we then have
\begin{align*}
\|\cm_\ct\phi\|^q_{p,q}&\le\bigg[\int^1_0\bigg(\int^1_0[\phi^\ast(tu)]^qdv(t)
\bigg)^{1/q}du\bigg]^q \\
&=\bigg[\int^1_0\bigg(\int^1_0[\phi^\ast(tu)]^qt^{q/p-1}dt\bigg)^{1/q}du\bigg]^q \\
&=\bigg[\int^1_0\bigg(\int^u_0[\phi^\ast(t)]^q\frac{t^{q/p-1}}{u^{q/p-1}}
\cdot\frac{dt}{u}\bigg)^{1/q}du\bigg]^q\\
&=\bigg[\int^1_0u^{-1/p}\bigg(\int^u_0[\phi^\ast(t)]^qt^{q/p-1}dt\bigg)^{1/q}du\bigg]^q \\
&\le\|\phi\|^q_{p,q}\bigg[\int^1_0u^{-1/p}du\bigg]^q=\bigg(\frac{p}{p-1}\bigg)^q\cdot\|\phi\|^q_{p,q},
\end{align*}
and so
\begin{eqnarray}
\|\cm_\ct\phi\|_{p,q}\le\frac{p}{p-1}\|\phi\|_{p,q}, \ \ \text{for} \ \ \phi\in L^{p,q}, \ \ q<p.  \label{eq4.8}
\end{eqnarray}
We end now the proof of Theorem \ref{thm4.3}.

Let $g:(0,1]\ra\R^+$ be non increasing. Then by Theorem \ref{thm3.1}
\[
\sup_{\phi^\ast=g}\|\cm_\ct\phi\|_{p,q}=\bigg[\int^1_0\bigg(\frac{1}{t}\int^t_0g(u)du\bigg)^q
dv(t)\bigg]^{1/q}
\]
so in order to prove that (\ref{eq4.8}) is sharp we just need to construct for every a such that $-\frac{1}{p}<a<0$, a non increasing $g_a:(0,1]\ra\R^+$ such that
\[
I/II\ra\bigg(\frac{p}{p-1}\bigg)^q, \ \ \text{as} \ \ a\ra-\frac{1^+}{p} \ \ \text{where}
\]
\[
I=\int^1_0\bigg(\frac{1}{t}\int^t_0g_a(u)du\bigg)^qt^{q/p-1}dt, \ \ \text{and}
\]
\[
II=\int^1_0[g_a(u)]^qt^{q/p-1}dt.
\]

But for $g_a(t)=t^a$, for a: $-\frac{1}{p}<a<0$ we have that
\[
I=\bigg(\frac{1}{a+1}\bigg)^q\frac{1}{q\big(a+\frac{1}{p}\big)} \ \ \text{and} \ \ II=\frac{1}{q\big(a+\frac{1}{p}\big)},
\]
so that
\[
I/II=\bigg(\frac{1}{a+1}\bigg)^q\overset{a\ra-\frac{1^+}{p}}{\longrightarrow}
\bigg(\frac{p}{p-1}\bigg)^q,
\]
and Theorem \ref{thm4.3} is proved.  \hfill$\square$
\end{Proof}
{\tiny{}}

\end{document}